\documentclass[12pt]{article}
\title{Compression of root systems and the $E$-sequence}
\author{Kevin Purbhoo\\ University of British Columbia}

\usepackage[mathscr]{eucal}
\usepackage{latexsym}
\usepackage{amssymb}
\usepackage{amsthm}
\usepackage{amscd}
\usepackage{amsmath}
\usepackage{epsfig}
\usepackage{graphics}

\newcommand{\ZZ}{\mathbb{Z}}
\newcommand{\RR}{\mathbb{R}}
\newcommand{\bolds}{\mathbf{s}}

\newcommand{\ahalf}{{\textstyle \frac{1}{2}}}
\newcommand{\link}{\mathcal{L}}
\newcommand{\antilink}{\mathcal{L}^c}
\newcommand{\veczero}{{\overrightarrow{0}}}
\newcommand{\veczeroset}{\{{\overrightarrow{0}}\}}

\newcommand{\orderideal}{\mathcal{J}}
\newcommand{\Dyn}{\mathsf{Dyn}}
\newcommand{\Dynhat}{\widehat{\mathsf{Dyn}}}
\newcommand{\Hasse}{\mathsf{H}}

\newtheorem{theorem}{Theorem}[section]
\newtheorem{lemma}[theorem]{Lemma}
\newtheorem{proposition}[theorem]{Proposition}
\newtheorem{corollary}[theorem]{Corollary}
\newtheorem{example}{Example}[section]
\newtheorem{definition}[example]{Definition}
\newtheorem{remark}[example]{Remark}
\newtheorem*{remark*}{Remark}
\newtheorem*{definition*}{Definition}

\begin{document}
\maketitle

\begin{abstract}
We examine certain maps from root systems to 
vector spaces over finite fields.  
By choosing appropriate bases, the images of these
maps can turn out to have nice combinatorial properties, which
reflect the structure of the underlying root system.  The main
examples are $E_6$ and $E_7$.
\end{abstract}


\section{Introduction}

The primary goal of this paper is to provide a convenient
way of visualising the root systems $E_6$ and $E_7$.  There
are two important relations on a root system that one might
wish to have a good understanding of: the poset structure,
in which $\alpha > \beta$ if $\alpha-\beta$ is a positive root,
and the orthogonality structure, in which $\alpha \sim \beta$
if $\alpha$ and $\beta$ are orthogonal roots. 

\begin{figure}[htbp]
  \begin{center}
    \epsfig{file=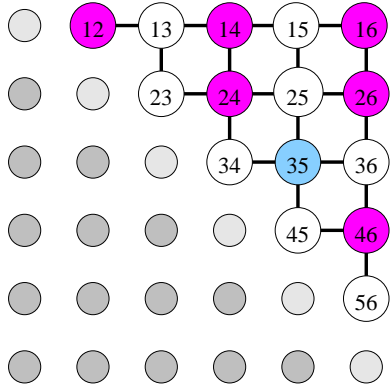, height=2in}
    \caption{Orthogonality and partial order in type $A$.}
    \label{fig:staircase}
  \end{center}
\end{figure}

In our paper on cominuscule Schubert calculus, with Frank Sottile 
\cite{PS}, we found that 
our examples required a good simultaneous understanding
both these structures.
This is easy enough to acquire for
the root systems corresponding to the classical Lie groups.  
In $A_n$, for example, one
can visualise the positive roots as the entries of an
strictly upper triangular $(n+1) \times (n+1)$ matrix, where
the $ij$ position represents the root $x_i-x_j$.
Then $\alpha \geq \beta$, if and only if $\alpha$ is weakly
right and weakly above $\beta$.  Orthogononality is also
straightforward in this picture:  $\alpha$ and $\beta$
are non-orthogonal if there is some $i$ such that crossing
out the $i^\text{th}$ row and the $i^\text{th}$ column
succeeds in crossing out both $\alpha$ and $\beta$. 
See Figure \ref{fig:staircase}.

In type $E$, it is less obvious how to draw
such a concrete picture.  Separately the two structures have
been well studied in the contexts of minuscule posets 
\cite{Proctor, Stembridge}, and strongly regular 
graphs (see e.g. \cite{BCN, Cameron, GR}).  However, once
one draws the Hasse diagram of the posets, the orthogonality
structure suddenly becomes mysterious.  Of course, one can
always calculate which pairs of roots are orthogonal, but
we would prefer a picture which allows us to do it instantly.
Thus the main thrust of this paper is to get to 
Figures~\ref{fig:cubecorner} and~\ref{fig:square}, which
illustrate how one can simultaneously visualise $E_7$ and $E_6$
posets and orthogonality structures, at least restricted to 
certain strata of the root system.  The restriction of these 
structures to the strata is exactly what is needed for the 
type $E$ examples in \cite{PS}.  With a little more work, one
can use these figures to recover the partial order
and orthogonality structures for the complete root system.

To reach these diagrams, we begin by considering certain maps from 
a root system to $(\ZZ/p)^m$, which are injective (or 2:1 if $p=2$).
Once we have some general observations about these maps, we give
examples for $E_6$ and $E_7$ which
are particularly nice.  In these cases,
we show that properties of the underlying root system
are reflected in simple combinatorial structures on the
target space, which is what allows us to produce diagrams in
question.  As the $E_7$ example is richer, we will
discuss it before the $E_6$ example.

The idea of relating the $E_6$ and $E_7$ root systems to 
$(\ZZ/p)^m$ has appeared elsewhere.  For example, Harris \cite{Harris}
uses such an identification to describe the Galois group of the
27 lines on the cubic surface---one of the del Pezzo surfaces.
The connection between del Pezzo surfaces and the exceptional Lie groups 
has been well established; we refer the reader to \cite{Manin}.  One
can also see such a relationship reflected in the well known 
identification of Weyl groups (see e.g. \cite{Atlas}):
\begin{gather*}
W(E_6) \cong SO(5; \ZZ/3) \cong O^-(6; \ZZ/2)  \\
W(E_7) \cong \ZZ/2 \rtimes Sp(6; \ZZ/2) .
\end{gather*}
These facts follow from the identifications outlined in this
paper, and presumably have been proved in similar ways before.

\section{Compression of root systems}

\subsection{Simply laced root systems}
Let $\Delta \subset \RR^n$ be a simply laced root system, so that
with respect to the inner product $\langle \cdot, \cdot \rangle$
on $\RR^n$ we have
$\langle \beta, \beta \rangle = 2$ for all $\beta \in \Delta$.  
We assume that $\Delta$ has full rank in $\RR^n$.
Let $\Lambda= \ZZ \Delta$ denote the lattice in $\RR^n$ generated by 
$\Delta$.

Choose a basis of simple roots 
$\alpha_1, \ldots, \alpha_n \in \Delta$, for $\Lambda$.  Let
$\Delta^+$ denote the positive roots with respect to this
basis, and $\Delta^-$ denote the negative roots.  Recall that
$\Delta^+$ is a partially ordered set, with $\beta > \beta'$
iff $\beta - \beta' \in \Delta^+$.  Roots $\beta$ and $\beta'$
are comparable in the partial ordering iff 
$\langle \beta, \beta' \rangle >0$

For each $\beta \in \Lambda$, we define $\beta^i$ to be the
coefficient of $\alpha_i$, when $\beta$ is expressed in the basis
of the simple roots: $\beta = \sum_{i=1}^n \beta^i\alpha_i$.

Let $\Dyn$ denote the Dynkin diagram Dynkin diagram of $\Delta$.  
As $\Delta$ is simply
laced, each component of $\Dyn$ has type ADE.  The vertices of 
$\Dyn$ are denoted $v_1, \ldots, v_n$, and correspond (respectively)
to the simple roots $\alpha_1, \ldots, \alpha_n$.  When
$\Delta$ is a simple root system (i.e $\Dyn$ has just one component),
the affine Dynkin diagram $\Dynhat$ is obtained by adding a vertex
$\hat v_n$ to $\Dyn$, corresponding to the lowest root $\hat \alpha_n$
of $\Delta$.

\subsection{Root systems over $\ZZ/p$}
Let $p\geq 2$ be a positive integer.  For reasons explained later
in this section,
we shall be mostly interested
in the case where $p$ is a prime, or $p=4$.
Let $V$ be a finite rank free module over $\ZZ/p$, with
a symmetric bilinear form $(\cdot|\cdot)$ taking values
in $\ZZ/p$.  Let 
$\Gamma = \{x \in V\setminus \veczeroset\ |\ (x|x)=2\}.$

Suppose that $\Gamma$ has 
subset $S = \{s_1, \ldots, s_n\}$ 
such that 
\begin{equation}\label{eqn:samedynkin}
\begin{alignedat}{2}
& &\qquad (s_i|s_j) &= \langle \alpha_i, \alpha_j \rangle \pmod p,
\quad \text{for all $i,j$,} \\
\text{and if $p=2$, assume} && \qquad s_i  &\neq s_j,
 \quad  \text{for all $i \neq j$.}
\end{alignedat}
\end{equation}
Then we obtain a map $f: \Lambda \to V$ by extending the
natural map $\alpha_i \mapsto s_i$ to a homomorphism of Abelian groups.

\begin{proposition}
If $\beta, \beta' \in \Lambda$ then 
\begin{equation}
(f(\beta)|f(\beta')) =
\langle \beta, \beta' \rangle \pmod p
\end{equation}
\end{proposition}

\begin{proof}
This is true for all pairs of simple roots, and both inner products 
are bilinear.
\end{proof}

\begin{corollary}
Suppose $\beta \neq \beta' \in \Delta$.  Then
$\langle \beta, \beta' \rangle = 0$ $\iff$
$(f(\beta)| f(\beta')) = 0$.
\end{corollary}

\begin{proof}
Since $\beta, \beta'$ are roots of a simply laced root system, 
$\langle \beta, \beta' \rangle \in \{-1,0,1\}$, thus 
$\langle \beta, \beta' \rangle = 0$ $\iff$
$\langle \beta, \beta' \rangle = 0 \pmod p$
$\iff$ $(f(\beta)|f(\beta')) = 0$.
\end{proof}

We now restrict the domain of $f$ to $\Delta$ if $p > 2$ and
$\Delta^+$ if $p=2$.

\begin{theorem}\label{thm:injective}
If $p >2$,
the map
$f: \Delta \to V$ is injective, and its image lies 
in $\Gamma$.
If $p =2$, the map
$f: \Delta^+ \to V$ is injective, and its image lies 
in $\Gamma$.
\end{theorem}

\begin{proof}
We first suppose $p>2$.
Note that the fact that
$f(\Delta) \subset \Gamma$ is clear from the fact that
every $\beta \in \Delta$ satisfies $\langle \beta, \beta \rangle = 2$.

Now, suppose that $\beta, \beta' \in \Delta$,
$f(\beta) = f(\beta')$.  We show that $\beta = \beta'$.

For all $\gamma \in \Delta$ we have 
$(f(\beta)| f(\gamma)) = (f(\beta')| f(\gamma))$, so 
$\langle \beta, \gamma \rangle = \langle \beta', \gamma \rangle
\pmod p$.  In particular the set of roots perpendicular
to $\beta$ and $\beta'$ are equal.
 implies $\beta$ and
$\beta'$ belong to the same simple component of $\Delta$.

There are two cases: if the component is of type $A_2$, then 
it is trivial that we must have $\beta = \beta'$.  If the component is 
not of type $A_2$, then the fact that $\beta$ and $\beta'$ have
the same set of perpendicular roots 
implies that $\beta = \pm \beta'$.  
(In types $D$ and $E$, the roots perpendicular 
to any given root span an entire hyperplane,
and in type $A$ it is easily checked.)
However, for all $x \in \Gamma$, $x \neq -x$.  
Since $f(\beta) = f(\beta') \in \Gamma$, we cannot have
$\beta' = -\beta$.  Thus $\beta = \beta'$.

For $p=2$, the fact that every $\beta \in \Delta^+$ satisfies 
$\langle \beta, \beta \rangle = 2$, implies that 
$f(\Delta^+) \subset \Gamma \cup \veczeroset$.  
It is therefore enough to show that
$f:\Delta^+ \cup \veczeroset \to \Gamma \cup \veczeroset$ is
injective.

Suppose $\beta, \beta' \in \Delta^+ \cup \veczeroset$,
$f(\beta) = f(\beta')$.  We show that $\beta = \beta'$.

As in the $p>2$ case, for all $\gamma \in \Delta$, we have
$\langle \beta, \gamma \rangle = \langle \beta', \gamma \rangle
\pmod 2$. Thus the sets $P(\beta)$ and $P(\beta')$, where
$$P(\beta):=\{\gamma \in \Delta\ |\ \langle \beta, \gamma \rangle =0\} 
\cup \{\pm \beta\},$$
coincide.  
This implies $\beta$ and
$\beta'$ belong to the same simple component of $\Delta$.
(In particular, if $f(\beta) = f(\beta') = \veczero$ then $\beta =\beta' =\veczero$.)

If this component is $A_1$ or $A_2$, it is trivial that 
$\beta = \beta'$.

If the component is $A_k$, $k \geq 4$, then $P(\beta)$ is a root 
system
of type $A_{k-2} \times A_1$,  where $\beta, \beta'$ are
both in the $A_1$ component.  
If the component is $E_6$, $E_7$, or $E_8$, then 
$P(\beta)$ is a root system of type $A_5 \times A_1$, $D_6 \times A_1$
or $E_7 \times A_1$,
where $\beta, \beta'$ are both in the $A_1$ component.
Thus in both these cases $\beta = \beta'$.

However, if the component if $D_k$, $k \geq 3$,
then $P(\beta)$ is a root
system of type $D_{k-2} \times A_1 \times A_1$, where $\beta,
\beta'$ are both in an $A_1$ component (a priori, not necessarily
the same one).  If $\beta, \beta'$ belong to the same
$A_1$ component, then $\beta = \beta'$.  So suppose they do not.
We identify the roots of $D_k$ with the vectors 
$\{e_i \pm e_j\ |\ i \neq j\} \subset \RR^k$, where $e_1, \ldots, e_k$
is an orthonormal basis for $\RR^k$.  It is easy to
see that if $\beta = e_i \pm e_j$, then $\beta' = e_i \mp e_j$.
So $\veczero= f(\beta) - f(\beta') = f(2e_j)$.  On the other hand,
for all $k$, we have $\veczero = 2f(e_k - e_j) + f(2e_j) = f(2e_k)$.
So $f(e_k \pm e_l) = f(e_k \mp e_l)$, for all $k \neq l$. But
among these must be a pair of simple roots.  We conclude that
$f$ restricted to the simple roots is not injective, hence $S$
contains fewer than $n$ elements, a contradiction.

\end{proof}

\begin{remark} \rm
Although we will not have use for it here, if $p$ is not a prime, one
could also allow the possibility that $V$ is not a free module.
In this case Theorem \ref{thm:injective} remains true provided
$f(\beta) \neq f(-\beta)$ for all $\beta \in \Delta$.  This will
be the case whenever $2 \nmid p$ or when $\Dyn$ has no 
component of type $A_1$.
\end{remark}

We now show that the most interesting cases
are when $p$ is a prime or $p=4$.
Suppose $p \geq 6$ is not prime.  Let $p' \neq 2$ be
a proper divisor of $p$.  Let $V' = V \otimes_{\ZZ/p} \ZZ/p'$.
Let $\rho:V \to V'$ denote the reduction modulo $p'$ map.
$V'$ comes with a symmetric bilinear form $(\cdot|\cdot)'$, the reduction of 
$(\cdot|\cdot)$ modulo $p'$.

\begin{corollary}
The composite map $f' := \rho \circ f : \Delta \to V'$ is injective,
and its image lies in $\Gamma' = \{x \in V'\ |\ (x|x)' = 2\}$.
Moreover $\langle \beta, \beta'\rangle = (f'(\beta)|f'(\beta'))' \pmod {p'}$.
\end{corollary}

\begin{proof}
As $p' \neq 2$, this follows from the fact that $\rho$ preserves inner
products modulo $p'$.
\end{proof}

\subsection{Compression}
The most interesting case of Theorem \ref{thm:injective} occurs
when rank $m$ of $V$ is smaller than the rank $n$ of 
$\Lambda$.  If this is the case, we will call the map $f$ a compression 
of the root system.  Here we give a necessary and nearly sufficient
condition for compression to be possible.

Let $A$ be the Coxeter matrix of $\Delta$,
$A_{ij} = \langle \alpha_i, \alpha_j \rangle$.

\begin{proposition}\label{prop:pdivides}
If we have $S$ as in equation \eqref{eqn:samedynkin}, and
$m < n$, then $p$ divides $\det(A)$.
\end{proposition}

\begin{proof}
Let $\bolds$ be the $m \times n$ matrix whose columns are are 
$s_i$ in some basis, and let $g$ be the $m \times m$ matrix
representing the bilinear form $(\cdot|\cdot)$ in the same
basis.  Then
$$A_{ij} = (s_i|s_j) = (\bolds^Tg\bolds)_{ij} \pmod p.$$
If $m<n$ then $\det(\bolds^Tg\bolds)= 0$, so $p|\det(A)$.
\end{proof}

Conversely, if $p$ is prime and $p|\det(A)$, and $A_p$ denotes the 
reduction of $A$ modulo $p$, then one can define $V = (\ZZ/p)^n/\ker(A_p)$.
and $s_i$ is the image of the standard basis vector $e_i$
under the natural map.  This will satisfy \eqref{eqn:samedynkin},
provided the $s_i$ are all distinct and non-zero.  The same construction
works if $p$ is not prime, though $V$ will not necessarily be
a free $\ZZ/p$-module.

In particular, we cannot hope for compression
in $E_8$, a root system for which $\det(A) = 1$.
For $E_7$, however, $\det(A) = 2$, and for $E_6$, $\det(A) = 3$.
Thus we should expect compression of the $E_7$ and
$E_6$ root systems to be possible, taking $p=2$ or $3$ 
respectively.

\subsection{Structures on $V$}
\begin{definition} \rm
The {\bf O-graph} of $V$ is the graph whose vertex set is $V$ and
whose edges are pairs $(x,y)$, $x \neq y$ such that $(x|y)=0$. 
The {\bf N-graph} of $V$ is the complement of the O-graph,
having vertex set $V$ and
edges $(x,y)$ such that $(x|y) \neq 0$. 
\end{definition}

As our two main examples involve $p=2$ and $p=3$, we consider
some special inner products $(\cdot|\cdot)$ in these case.

If $p=2$, we
let $V$ be an even dimensional vector space over $\ZZ/2$ with a 
symplectic form $(\cdot|\cdot)$.  By symplectic form, we mean an
(anti)symmetric non-degenerate bilinear form for which $(x|x) = 0$ for 
all $x$.
Thus $\Gamma = V \setminus \veczeroset$.
 We see that $S \subset V \setminus \veczeroset$ satisfies the condition
\eqref{eqn:samedynkin} iff the restriction of the N-graph to
$S$ is isomorphic to $\Dyn$.  In this case, the associated
map $f$ gives an injective map from $\Delta^+$ to 
$V \setminus \veczeroset$.

If $p=3$, we take $V$ to be an $m$-dimensional vector space over 
$\ZZ/3$, with the standard symmetric form 
\begin{equation}\label{eqn:symmform3}
\big((x_1, \dots, x_m)\big|(y_1, \ldots, y_m)\big) = \sum_{i=1}^m x_iy_i.
\end{equation}
Note that
\begin{equation}
\Gamma 
= \big\{(x_1, \ldots, x_m)
\ \big|\ \#\{i\;|\;x_i \neq 0\} = 2 \pmod 3 \big\}.
\end{equation}

\section{Application to type $E$}

\subsection{The $E$-sequence}
Consider $\RR^8$ with the standard Euclidean inner product
$\langle \cdot, \cdot \rangle$.  Let
$\alpha_1, \ldots, \alpha_8$ be the vectors
\begin{equation*}
\begin{aligned}[c]
\alpha_1 &= (1,-1,0,0,0,0,0,0) \\
\alpha_2 &= (\ahalf,\ahalf,\ahalf,-\ahalf,-\ahalf,-\ahalf,-\ahalf,-\ahalf) \\
\alpha_3 &= (0,1,-1,0,0,0,0,0) \\
\alpha_4 &= (0,0,1,-1,0,0,0,0) \\
\alpha_5 &= (0,0,0,1,-1,0,0,0) \\
\alpha_6 &= (0,0,0,0,1,-1,0,0) \\
\alpha_7 &= (0,0,0,0,0,1,-1,0) \\
\alpha_8 &= (0,0,0,0,0,0,1,-1).
\end{aligned}
\qquad \qquad
\epsfig{file=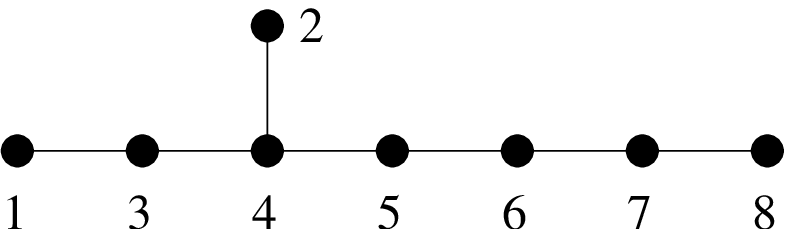, height=.6in}
\end{equation*}
These are the simple roots of $E_8$, which span the $E_8$ lattice.  
They correspond to the vertices $v_1, \ldots, v_8$ of $\Dyn(E_8)$, in 
the order shown on the right.
%
%
This ordering of simple roots of $E_8$ corresponds
to the inclusion of groups below.
$$
\begin{CD}
A_1  @>>> D_2 @>>> E_3 @>>> E_4 @>>> E_5 @>>> E_6 @>>> E_7 @>>> E_8 \\
@.        @.       @|        @|       @|       @.      @.      @.  \\
 @.    @. A_2 \times A_1 @>>> A_4 @>>> D_5 @. @. @. \\
\end{CD}
$$

To obtain the root systems
of $E_n$, $3 \leq n \leq 7$ we take the simple roots to be
$\alpha_1, \ldots, \alpha_n$. These span the $E_n$ lattice.  In
general, the roots of $E_n$ are the Lattice vectors 
$\alpha \in \Lambda$ such that $\langle \alpha, \alpha \rangle = 2$.

The roots $\Delta(E_8)$ of $E_8$, are stratified as 
$\Delta(E_8) = \coprod \Delta_s$.
For $s \neq 2,3$,
\begin{equation}\label{eqn:strata}
\Delta_s = \{\beta \in \Delta(E_8)\ |\ 
\text{$\beta \geq \alpha_s$, and  $\beta \ngeq \alpha_t$ for all $t >s$} \}.
\end{equation}
$\Delta_s$ are the roots of $E_s$ minus the roots
of $E_{s-1}$.  Equation \eqref{eqn:strata} makes sense for $s=2$, and $s=3$;
however it is convenient (and arguably correct) to put 
these into the same stratum:
$$\Delta_3 = \{\pm \alpha_3, \pm \alpha_2, \pm(\alpha_1+ \alpha_2)\}.$$
We also have a stratification of $\Delta^+(E_8) =\coprod \Delta_s^+$ given
by $\Delta_s^+ = \Delta^+(E_8) \cap \Delta_s$.
This stratification has the property that the roots of $E_n$,
$3 \leq n \leq 8$ are precisely 
$$\Delta(E_n) = \coprod_{s\leq n} \Delta_s.$$

For notational convenience, we define $s' = \max\{3, s+1\}$,
so that $\Delta_s$ and $\Delta_{s'}$ always denote consecutive strata.

\begin{figure}[htbp]
  \begin{center}
    \epsfig{file=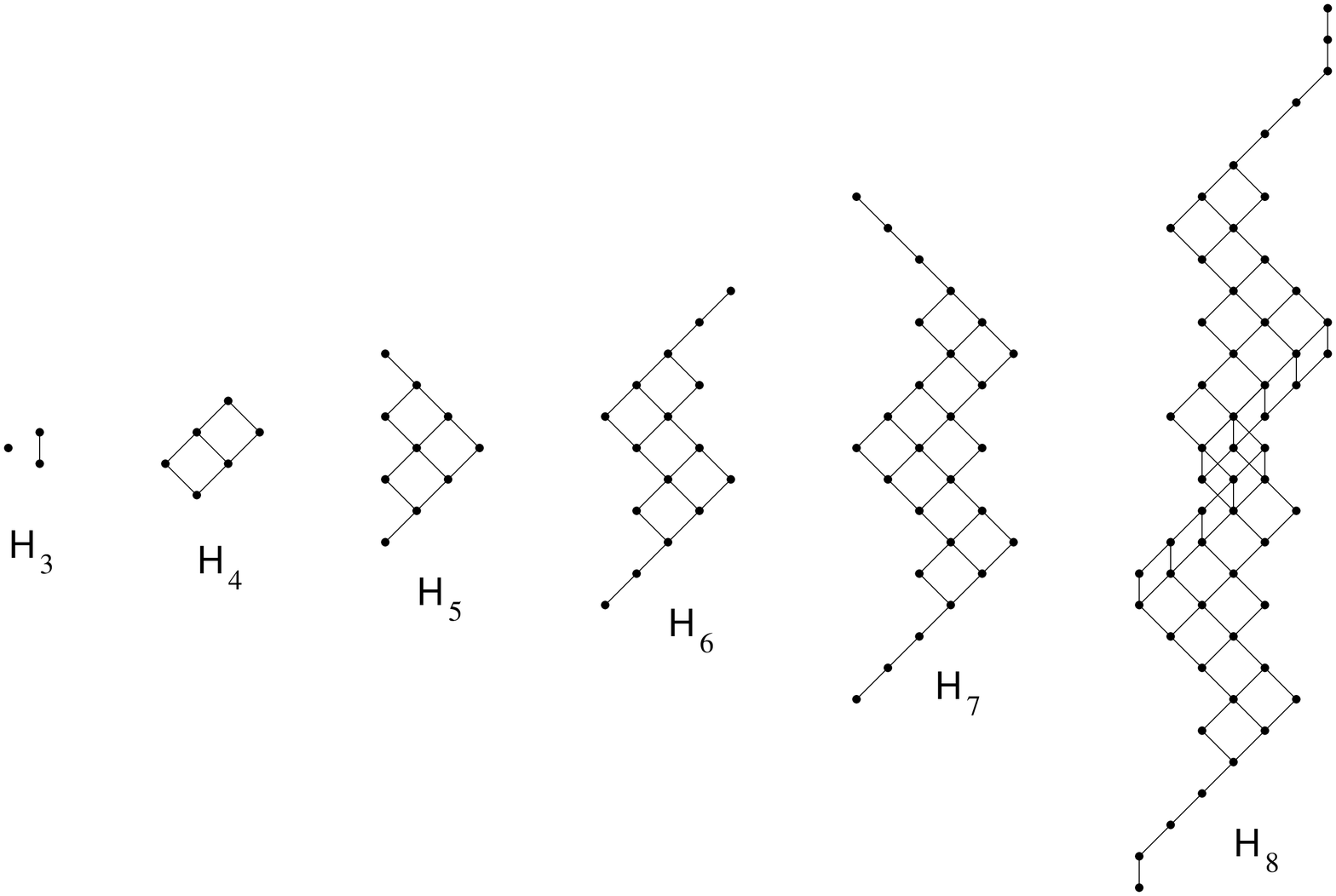, height=3.25in}
    \caption{The Hasse diagrams $\Hasse_s$, $3 \leq s \leq 8$.}
    \label{fig:hasse}
  \end{center}
\end{figure}

For each stratum let $\Hasse_s$ denote the graph whose vertices are
$\Delta_s^+$ and whose edges form the Hasse diagram of the poset
structure on $\Delta^+$, restricted $\Delta_s^+$.  Thus
we have an edge joining $\beta$ and $\beta'$ if one of 
$\beta \pm \beta'$ is a simple root.  These are shown in Figure
\ref{fig:hasse}.

Finally, it is worth noting the size of each stratum.  The
stratification $\Delta^+(E_8) =\coprod \Delta_s^+$  has strata
of sizes 
$1\ (s=1)$, $3\ (s=3)$, $6\ (s=4)$, $10\ (s=5)$, $16\ (s=6)$, $27\ (s=7)$ 
and $57\ (s=8)$.

\subsection{A compression of $E_7$}
We now take $\Delta$ to be the $E_7$ root system.

Let $F = (\ZZ/2 \times \ZZ/2, \oplus)$ denote the non-cyclic four element 
group.  
We denote the elements of this group $\{0,1,2,3\}$, and the operation
$a\oplus b$ is binary addition without carry.
$F$ is a two dimensional vector space over $\ZZ/2$ and thus has
a unique symplectic form:
$$(a|a') = \begin{cases}
0 &\text{if $a=0$, $a'=0$ or $a=a'$} \\
1 &\text{otherwise}.
\end{cases}
$$

We shall take $V= F^3$, and whenever possible
we write a triple $(a,b,c) \in V$ simply as $abc$.   We endow $V$
with the symplectic form
$$(abc|a'b'c') = (a|a') + (b|b') + (c|c').$$

We take as our subset $S \subset V$, the set $S = \{s_1, \ldots,
s_7\}$, where
\begin{equation*}
\begin{gathered}[b]
\begin{alignedat}{2}
s_1 &= 100, &\quad s_2 &= 030, \\
s_3 &= 300, &\quad s_4 &= 111,
\end{alignedat} \\
s_5 = 003, \quad s_6 = 001, \quad s_7 = 033.
\end{gathered}
\qquad \qquad
\epsfig{file=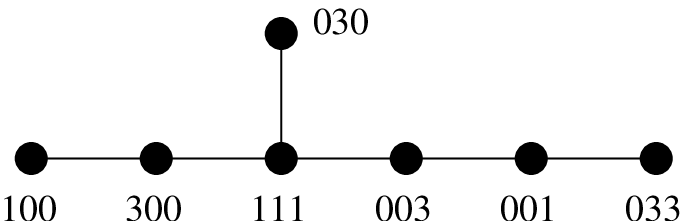, height=0.6in}
\end{equation*}

\begin{proposition}
The restriction of the N-graph to $S$ is $\Dyn(E_7)$.
The natural homomorphism $f$ takes $\alpha_i$ to $s_i$.
\end{proposition}

\begin{proof}
This just needs to be checked.
\end{proof}

As a consequence of we obtain the following corollary of Theorem
\ref{thm:injective}.

\begin{corollary}
The map $f: \Delta^+ \cup \veczeroset \to V$ is a bijection.
\end{corollary}

\begin{proof}
It is an injection by Theorem \ref{thm:injective}.  But
$\#(\Delta^+ \cup \veczeroset) = \#(V) = 64$, thus it is a bijection.
\end{proof}

\subsection{Restriction to strata}
\label{sec:strata}

Let $\Gamma_s^+$ denote the image of the stratum $\Delta_s^+$ under
$f$.  Here we show how natural structures on $\Delta_s^+$
are preserved under $f$, and are more palatable in $\Gamma_s^+$.

We define a new graph structure on $V$.
\begin{definition} \rm
The {\bf T-graph} is the graph with vertex set 
$V = F^3$, and $abc$ adjacent to $a'b'c'$, if exactly one of 
$\{a=a',b=b',c=c'\}$
holds.  
\end{definition}

Unlike the O-graph, the the T-graph has translation symmetries:
for any $x \in V$, the map 
\begin{equation}\label{eqn:aut1}
y \mapsto x \oplus y 
\end{equation}
is an automorphism.  It is a strongly regular graph.  In particular 
every vertex has valence 27.

\begin{definition} \rm
For $v \in V$, the {\bf link} on $v$ in the T-graph, denoted
$\link(v)$, is the
restriction of the T-graph to the vertices adjacent to $v$.
The {\bf antilink} on $v$ in the T-graph, denoted $\antilink(v)$
is the restriction of the T-graph of the vertices non-adjacent
to $v$ (not including $v$).
\end{definition}

\begin{lemma}\label{lem:27}
The image of the largest stratum $\Gamma_7^+$ is the vertex set
of $\link(\veczero)$.  
\end{lemma}
In other words,
$\Gamma_7^+$ is the set of $abc \in V$ such that exactly one of
$\{a=0,b=0,c=0\}$ holds.

Note that this result is not independent of the choice of $S$
for the images of the simple roots.  We have chosen $S$ quite
carefully, in part to make this lemma hold.  It is possible,
(and not difficult) to check this result on each of the 27 roots 
of $\Delta_7^+$; however, since a symmetry argument is available, 
we present it here.

\begin{figure}[htbp]
  \begin{center}
    \epsfig{file=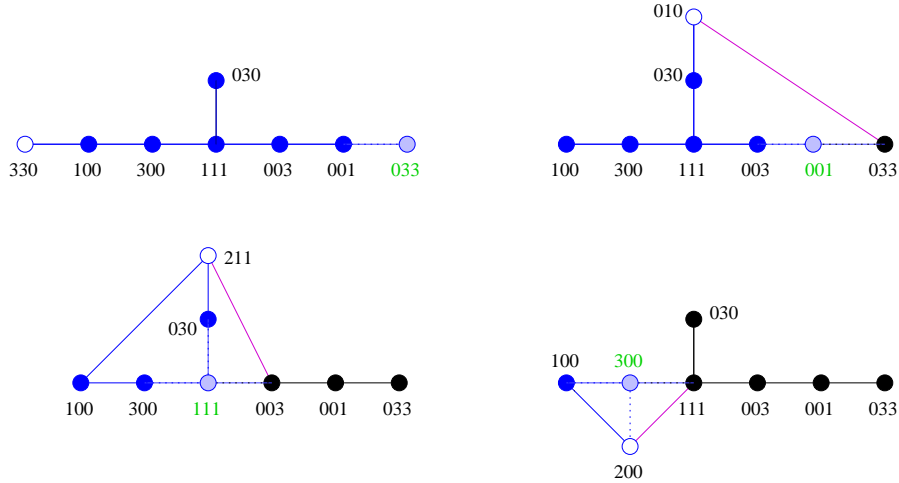, height=2.5in}
    \caption{The automorphisms of $\Gamma_7^+$ described in the
     proof of Lemma \ref{lem:27} in terms of the Dynkin diagram
     of $E_7$. }
    \label{fig:automorphisms}
  \end{center}
\end{figure}

\begin{proof}
We know that $f(\alpha_7) = 033 \in \Gamma_7^+$, thus it
suffices to show that $\Gamma_7^+$ is invariant under the
following symmetries:
\begin{equation}
\label{eqn:symmetries27}
\begin{aligned}
(a,b,c) &\mapsto (c,b,a) \\
(a,b,c) &\mapsto (a,c,b) \\
(a,b,c) &\mapsto ([1 \to 2 \to 3 \to 1]\cdot a,b,c) \\
(a,b,c) &\mapsto ([3 \leftrightarrow 2] \cdot a,b,c) 
\end{aligned}
\end{equation}

These symmetries come from a Dynkin diagram construction, which
we first describe for any $E_n$.  A similar construction can
also be used for types $A$ and $D$.
Let $D= \Dyn(E_n)$.  We decorate each vertex 
of $D$ with the corresponding simple root in
$\Delta$.

Choose a vertex $v_i \in D$, where $i \notin \{1,2,8\}$.  If we
delete the edge $(v_i, v_{i+1})$ from $D$, the diagram 
breaks up into two 
components $D', D''$ where $D'$ is the component containing
$v_1$.  If $i=n$, $D''$ will be empty.  Note that $D'$ is
a sub-Dynkin diagram of $D$, and
hence corresponds to a sub-root system $\Delta' \subset \Delta$.  

We apply the following construction to obtain a new Dynkin diagram 
$\tilde D$:
\begin{enumerate}
\item Add to $D'$ the affine vertex $\hat v$, to form the affine
Dynkin diagram $\widehat{D'}$.  This vertex
is decorated with the {\em lowest root} 
$\hat \alpha \in \Delta'$.
\item For every vertex $\widehat{D'}$, replace the root 
 which decorates the vertex by its negative.
\item Delete the vertex $v_i$.
\item If $D''$ is not empty, reattach it by forming an edge
 $(\hat v, v_{i+1})$.  The result is $\tilde D$.
\end{enumerate}

The underlying graph $\tilde D$ is isomorphic to $D$, but
under this isomorphism, the roots decorating the vertices have 
changed.  The roots decorating $\tilde D$  correspond to a new system 
of simple roots for $\Delta$.  Thus this process corresponds gives an 
automorphism of $\Delta$, and hence to an automorphism of $\Gamma$.

Returning now to the $E_7$ case, we note that for any root 
$\beta$, $\beta^7$ is preserved modulo
$2$ under each of these automorphisms.  Thus each 
automorphism restricts to an automorphism of $\Delta_7$,
and hence of $\Gamma_7^+ = f(\Delta_7^+) = f(\Delta_7)$.
The symmetries \eqref{eqn:symmetries27}, are the automorphisms of $\Gamma_7^+$
given by the construction above, using vertices  $v_7, v_6, v_4, v_3$ 
respectively.  It is sufficient
to check this on the images of the simple roots, and this
is easily done.  See Figure \ref{fig:automorphisms}.
\end{proof}

%
%
%
%
%
%
%
%
%
%
%

The following construction is useful for relating the other
strata to $\Delta_7^+$.
Put $z_7 = \veczero$, $\zeta_s = \sum_{i=s'}^7 \alpha_i$, and
$z_s = f(\zeta_s)$
for $s=1,3,4,5,6$.  If $\beta \in \Delta_s^+$, define
$\tilde \beta = \beta + \zeta_s$.  Note that 
$\tilde \beta \in \Delta_7^+$.

\begin{theorem}\label{thm:T-graph7}
On $\Gamma_s^+$, the image of any stratum, the T-graph agrees with
the O-graph.
In particular, for $\beta, \beta' \in \Delta_s^+$, 
$\langle \beta , \beta' \rangle = 0$ if and only if
$f(\beta)$ and $f(\beta')$ agree in exactly one coordinate.
\end{theorem}

\begin{proof}
We assume $\beta \neq \beta'$, since the result is trivial otherwise.
Put $f(\beta) = abc$, $f(\beta') = a'b'c'$.

We have
$
\langle \beta, \beta' \rangle = 
(a|a') + (b|b') + (c|c') \pmod 2
$
Thus $\langle \beta, \beta' \rangle = 0$ $\iff$ an odd number of
$\{(a|a'), (b|b'), (c|c')\}$ are non-zero.  

We first show the theorem is true for the stratum $\Delta_7^+$.
By Lemma
\ref{lem:27} exactly one of $\{a,b,c\}$ and exactly one of
$\{a',b',c'\}$ is zero.  Suppose $a=0, a'=0$. Then we have
$\langle \beta, \beta' \rangle = 0$ $\iff$ $b \neq b'$ and $c \neq c'$
$\iff$ $abc$ and $a'b'c'$ agree in exactly one coordinate.  
Suppose $a=0, b'=0$.  Then
$\langle \beta, \beta' \rangle = 0$ $\iff$ $c=c'$ $\iff$
$abc$ and $a'b'c'$ agree in exactly one coordinate.  
The remaining cases are the
same as these two by symmetry.

To show the statement for other strata $\Delta^+_s$, we 
consider $\tilde \beta$ and $\tilde \beta'$.
As $\beta, \beta'$ belong to the same stratum, we have
$\langle \tilde \beta, \tilde \beta' \rangle = 
\langle \beta, \beta' \rangle$.
Also,
$f(\tilde \beta) = f(\beta) \oplus z_s$
is the image of $f(\beta)$ under an automorphism of the T-graph;
thus $f(\tilde \beta)$ is adjacent to $f(\tilde \beta')$ $\iff$
$f(\beta)$ is adjacent to $f(\beta')$.

\end{proof}

In general, all of the strata can be described in terms
of links in the T-graph. 

\begin{theorem}\label{thm:otherstrata}
We have the following identifications:
\begin{enumerate} 
\item 
$\displaystyle 
\Gamma^+_s = \link(z_s) \cap \Big(\bigcap_{7 \geq t>s} \antilink(z_t) \Big ).$
\item 
$\displaystyle \Gamma^+_s \oplus z_s
= \link(z_7) \cap \Big(\bigcap_{7 > t \geq s} \antilink(z_t) \Big).$
\end{enumerate}
\end{theorem}

\begin{proof}
We already know this is true for $s=7$, so assume
$s \leq 6$.

Let us calculate $\link(z_t) \setminus \link(\veczero)$.
First note that by Theorem \ref{thm:T-graph7}
\begin{equation}
\label{eqn:7-t}
\link(\veczero) \setminus \link(z_t) = 
\{x \in \Gamma_7^+\ |\ (x|z_t) = 1\}
\end{equation}
We have 
$y \in \link(z_t) \setminus \link(\veczero)$
$\iff$ $y \oplus z_t \in \link(z_t \oplus z_t) \setminus \link(\veczero \oplus
z_t)
= \Gamma_7^+ \setminus \link(z_t)$
$\iff$ $y \oplus z_t \in \Gamma_7^+$ and $(y\oplus z_t|z_t) = 1$,
by \eqref{eqn:7-t}.
Now since $z_t \in \Gamma_7^+$, 
$(y|z_t) = (y\oplus z_t|z_t) = 1$ implies that
$y \oplus z_t \in \Gamma_7^+$ if and only if 
$y \notin \Gamma_7^+$.  Thus we see that
\begin{equation}
\label{eqn:t-7}
\link(z_t) \setminus \link(\veczero) =
\{y \in V\ |\ y \notin \Gamma_7^+,\ (y|z_t) =1\}.
\end{equation}

The set $\Gamma_s^+$ is characterized by 
$$y \in \Gamma_s^+\quad \iff \quad (y|z_i) =
\begin{cases}
1 & \text{if $i = s$} \\
0  & \text{if $i > s$.}
\end{cases}
$$
Thus, using \eqref{eqn:t-7} we see that
\begin{align*}
\Gamma_s^+ &=
\{y \in \Gamma\ |\ y \notin \Gamma_7,\ (y|z_s) =1,
\ (y|z_t)=0 \text{ for all $t>s$}
\} \\
&= \Big(\link(z_s) \setminus \link(\veczero) \Big)
\setminus \Big(
\veczeroset \cup \bigcup_{t=s+1}^6
\link(z_t) \setminus \link(\veczero) \Big) \\
&= \link(z_s) \cap \Big(\bigcap_{7 \geq t>s} \antilink(z_t) \Big ).
\end{align*}
On the other hand, using \eqref{eqn:7-t},
\begin{align*}
\Gamma_s^+ \oplus z_s 
&= \{x \in \Gamma_7^+\ |\ (x \oplus z_s|z_s) = 1,
\ (x \oplus z_s|z_t) = 0 \text{ for all $t>s$}
\} \\
&= \{x \in \Gamma_7^+\ |\ (x |z_t) = 1,
\text{ for all $t\geq s$} \} \\
&= \bigcap_{7 > t \geq s} \link(\veczero) \setminus \link(z_t)  \\
&= \link(z_7) \cap \Big(\bigcap_{7 > t \geq s} \antilink(z_t) \Big).
\end{align*}

\end{proof}

\subsection{Partial ordering}
\label{sec:order}

We now show how one can recover the partial ordering on $\Delta_s^+$
from $\Gamma_s^+$.

\begin{lemma}\label{lem:orthseq}
Let $x_1, x_2 \in \Gamma_7^+$.  If $x_1$ and $x_2$ are orthogonal, 
there exists a unique vector
$x_3 \in \Gamma_7^+$ such that $\{x_1, x_2,x_3\}$ are
pairwise orthogonal, and moreover, $x_1\oplus x_2\oplus x_3=\veczero$.
Conversely if $x_1 \oplus x_2 \in \Gamma_7^+$ then $x_1$ and $x_2$
are orthogonal.
\end{lemma}

In light of Theorem \ref{thm:T-graph7}, this is quite easy
to show for our preferred choice of $S$.  Nevertheless, this
result is true for any $S$ satisfying \eqref{eqn:samedynkin},
and so we give a more general proof.

\begin{proof}
Let $\beta_1 = f^{-1}(x_1) \in \Delta_7^+$
and $\beta_2 = f^{-1}(x_2) \in \Delta_7^+$. 
View $\beta_1$ and $\beta_2$ as roots in the $E_8$ root system.  
Assume that $x_1$ and $x_2$ are orthogonal; hence $\langle \beta_1,
\beta_2\rangle = 0$.
Throughout the proof we use the fact that the sum of two roots is 
a root if their inner product is negative.

To begin, for any $\beta \in \Delta_7^+$, we have 
$\langle \alpha_8, \beta \rangle = -1$
so $\alpha_8 + \beta$ is a root of $E_8$.
Similarly,
$\langle \alpha_8+ \beta_1, \beta_2 \rangle = -1$
so $\alpha_8 + \beta_1 + \beta_2$ is a root of $E_8$.

To show existence let $\hat \alpha_8$ denote the affine (lowest)
root of $E_8$.  Then 
$\langle \alpha_8+ \beta_1 + \beta_2, \alpha_8 + \hat \alpha_8 \rangle = -1$
so $\hat \alpha_8 + 2\alpha_8+\beta_1 + \beta_2$ is a root.
Let 
\begin{equation}
\label{eqn:defbeta3}
\beta_3 = - (\hat \alpha_8 + 2\alpha_8+\beta_1 + \beta_2),
\end{equation}
and $x_3 = f(\beta_3)$.
Note that $\beta_3^8=0$, $\beta_3^7 = 1$, so $\beta_3 \in \Delta_7^+$,
hence $x_3 \in \Gamma_7^+$.  And we can explicitly check
$\langle \beta_3, \beta_1\rangle = \langle \beta_3, \beta_2\rangle = 0$,
so $\{x_1, x_2, x_3\}$ are pairwise orthogonal.  
Finally, note that the affine root of $E_8$ has
the property that $\langle \hat \alpha_8, \gamma\rangle = 0$
for all roots $\gamma \in \Delta(E_7)$.  Thus, for all
$\gamma \in \Delta(E_7)$, we have
\begin{align*}
(x_1 \oplus x_2 \oplus x_3|f(\gamma))
&= \langle \beta_1 +\beta_2 + \beta_3, \gamma \rangle \pmod 2 \\
&= \langle \hat \alpha_8 - 2\alpha_8, \gamma \rangle \pmod 2 \\
&= 0 .
\end{align*}
Thus $x_1 \oplus x_2 \oplus x_3 = \veczero$.

For uniqueness, let $x_3$ be any vector orthogonal
to $x_1$ and $x_2$, and let $\beta_3= f^{-1}(x_3) \in \Delta_7^+$.
$\langle \alpha_8+ \beta_1 + \beta_2, \alpha_8 + \beta_3 \rangle = -1$,
so $\gamma = -(2\alpha_8 + \beta_1 + \beta_2 + \beta_3)$
is a root of $E_8$. But $(\gamma)^8 = -2$, and the only root of $E_8$
with this property is the affine root $\hat \alpha_8$.  We conclude
that \eqref{eqn:defbeta3} must hold.

For the converse, note that if $x_1$ and $x_2$ are not orthogonal,
then $\beta_1 - \beta_2$ is a root not in $\Delta_7$, hence
$x_1 \oplus x_2 \notin \Gamma_7^+$.
\end{proof}

\begin{corollary}\label{cor:orthseqcor}
Suppose $x,y \in \Gamma_s^+$. Then
$(x|y)=0$ if and only if  $x \oplus y \in \Gamma_7^+$
\end{corollary}

\begin{proof}
Note that
$x\oplus z_s, y\oplus z_s \in \Gamma_7^+$.  Thus
$(x|y)=0$ $\iff$ $(x\oplus z_s|y \oplus z_s) = 0$ 
$\iff$
$x\oplus y = x \oplus z_s \oplus y \oplus z_s \in \Gamma_7^+$, 
by Lemma \ref{lem:orthseq}.
\end{proof}

\begin{theorem}\label{thm:order7}
Let $\alpha \notin \Delta_7^+$ be a positive root.  Let 
$\beta \in \Delta_s^+$.  Then
$f(\beta)\oplus f(\alpha) \in \Gamma_s^+$ if and only if
either $\beta + \alpha \in \Delta_s^+$ or 
$\beta - \alpha \in \Delta_s^+$.
\end{theorem}

\begin{proof}
Certainly if one of $\beta \pm \alpha \in \Delta_s^+$, then
$f(\beta)\oplus f(\alpha) \in \Gamma_s^+$.  Suppose that
neither $\beta + \alpha$ nor $\beta -\alpha$ is in $\Delta_s^+$.
Note that neither can be in $\Delta_s^-$ either.
If one of $\beta \pm \alpha$ is a root, then it belongs to
some stratum other than $\Delta_s^+$, so 
$f(\beta)\oplus f(\alpha) \notin \Gamma_s^+$.
On the other hand if neither is a root then 
$\langle \beta, \alpha \rangle = 0$.  Suppose that
$f(\beta)\oplus f(\alpha) \in \Gamma_s^+$.  Then
\begin{align*}
(f(\beta)\oplus f(\alpha)|f(\beta)) &= (f(\alpha)|f(\beta)) \\
&= \langle \beta, \alpha \rangle \pmod 2 \\
&= 0.
\end{align*}
By Corollary \ref{cor:orthseqcor}, we conclude 
$f(\alpha) \in \Gamma_7^+$,
a contradiction.
\end{proof}

\begin{figure}[htbp]
  \begin{center}
    \epsfig{file=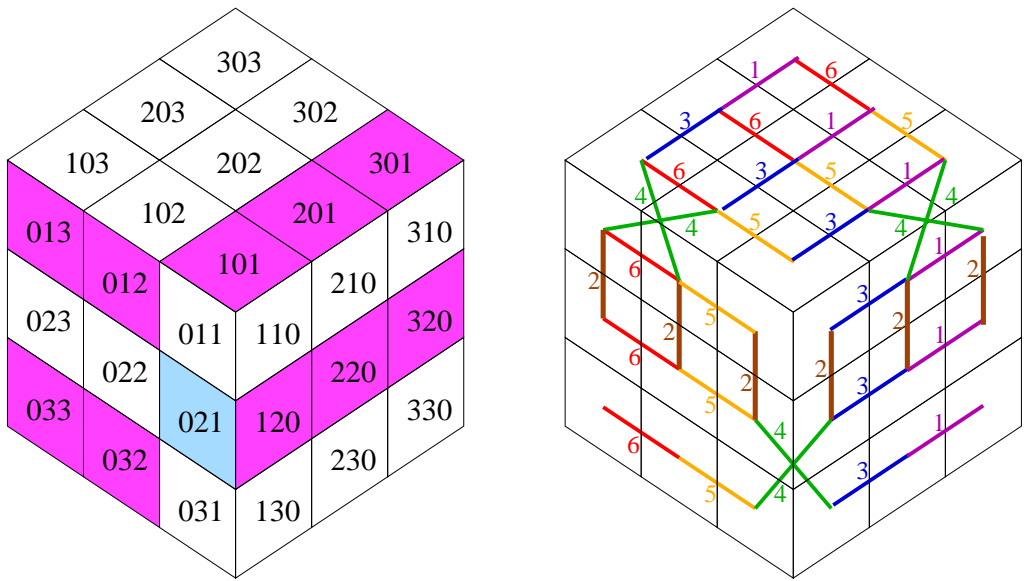, height=3in}
    \caption{Orthogonality and partial order in $\Gamma_7^+$. }
    \label{fig:cubecorner}
  \end{center}
\end{figure}

One can visualise $\Gamma_7^+$ as the squares one sees looking
at the
corner of a $3 \times 3 \times 3$ cube.  The elements are arranged
as shown in Figure \ref{fig:cubecorner}.  We can recover the 
Hasse diagram $\Hasse_7$ of the
poset structure on $\Delta_7^+$ in this picture, as follows.  
For each simple root $\alpha_i$, $i=1, \ldots, 6$, we draw an 
edge joining $x$ and $y$ if $y = x \oplus f(\alpha_i)$.
By Theorem~\ref{thm:order7} we will draw such an edge if
and only if the corresponding roots in $\Gamma_7^+$ are related by 
addition a simple root, which is exactly how the Hasse diagram
is constructed.  A similar procedure also works on the smaller
strata.

In this picture, orthogonality is easy to determine as well.
By Theorem \ref{thm:T-graph7} this is determined by 
the links in the T-graph.
For any $x \in \Gamma_7^+$, the set of $y \in \Gamma_7^+$ orthogonal 
to $x$ can be described as follows:  if $y$ is on the same face
of the $3 \times 3 \times 3$ cube as $x$, then $y$ is not in the same
row or column as $x$; if $y$ is on a different face from $x$
then $y$ is in the same extended row/column as $x$.
Figure \ref{fig:cubecorner} shows $\link(021) \cap \Gamma_7^+$,
which is set of root images orthogonal to $021$.

More generally, one can compare 
$\alpha \in \Delta_s^+$ and 
$\beta \in \Delta_t^+$ by considering $\tilde \alpha$, and 
$\tilde \beta$.

\begin{proposition} Suppose $\alpha \in \Delta_s^+$ and 
$\beta \in \Delta_t^+$, and $s \leq t$.  Then $\alpha < \beta$
if and only if $\tilde \alpha < \tilde \beta$.  
If $s = t$ then
$\alpha$ is orthogonal to $\beta$ if and only if 
$\tilde \alpha$ is orthogonal to $\tilde \beta$.
If $s<t$ then 
$\alpha$ is orthogonal to $\beta$ if and only if 
$\tilde \beta$ is orthogonal to both
$\tilde \alpha$ and $\zeta_s$, or to neither.
\end{proposition}

\begin{proof}
For the first statement, the `if' direction is clear,
as $\tilde \beta - \tilde \alpha > \beta - \alpha$.
Conversely, if $\beta> \alpha$ then $\beta^i \geq 1$ for
$i = s', \ldots, t$, hence 
$\tilde \beta - \tilde \alpha = \beta - \alpha - 
\sum_{i=s'}^t \alpha_i$ is still positive.

The statements about orthogonality follow from the
following calculation:
$$(f(\alpha)|f(\beta)) =
\begin{cases}
(f(\tilde \alpha)|f(\tilde \beta))
&\text{if $s=t$} \\
(f(\tilde \alpha)|f(\tilde \beta)) \oplus
(z_s|f(\tilde \beta)) &\text{if $s < t$}.
\end{cases}
$$
\end{proof}

\subsection{Action of the Weyl group of $E_6$ on $\Gamma_7^+$}

The O-graph restricted to $\Gamma_7^+$ is a well known object; its
complement is the Schl\"afli graph (see e.g. \cite{BCN, Cameron, GR} 
for alternate desciptions), which describes the
incidence relations of the 27 lines on a cubic surface.  It is well 
known that
the full automorphism group of 
the Schl\"afli graph is the Weyl group of $E_6$.  Many of these 
automorphisms are manifest from our description.

If $\phi_1, \phi_2, \phi_3$
are automorphisms of $F$, then 
\begin{equation}\label{eqn:aut2}
(a_1, a_2, a_3) \mapsto (\phi_1(a_1), \phi_2(a_2), \phi_3(a_3))
\end{equation}
is manifestly an automorphism of the Schl\"afli graph.  
If $\pi \in S_3$ is a permutation of 
$\{1,2,3\}$ then we have the automorphism
\begin{equation}\label{eqn:aut3}
(a_1, a_2, a_3) \mapsto (a_{\pi(1)}, a_{\pi(2)}, a_{\pi(3)}).
\end{equation}

If $\alpha \in \Delta(E_6)$,
then the action of the reflection $r_\alpha$ on $\Gamma_7^+$ is
given by 
$$r_\alpha(f(\beta)) = f(r_\alpha(\beta)) = 
\begin{cases}
\beta \pm \alpha &\text{if $\langle \alpha, \beta \rangle = \mp 1$} \\
\beta &\text{if $\langle \alpha, \beta \rangle = 0$}.
\end{cases}
$$
Using Theorem \ref{thm:order7}, we see that
for $x \in \Gamma_7^+$,
\begin{equation}\label{eqn:reflroot}
r_\alpha(x) =
\begin{cases}
x \oplus f(\alpha) &\text{if $x \oplus f(\alpha) \in \Gamma_7^+$} \\
x &\text{otherwise}.
\end{cases}
\end{equation}
Each $r_\alpha$ swaps six pairs $x \leftrightarrow x \oplus f(\alpha)$
and the restriction of the O-graph to these 12 vertices is a union
of two $K_6$ graphs, which are maximal cliques.  These pairs are
known as Schl\"afli double sixes---there are 36 in total, each
arising in this way for some unique $\alpha \in \Delta^+(E_6)$.

From \eqref{eqn:reflroot}, it is easy to verify that the automorphisms 
\eqref{eqn:aut2} are generated by reflections in the roots 
$\alpha_1, \alpha_3$ (generating all possible $\phi_1$);  
$\hat \alpha_6, \alpha_2$ (for $\phi_2$);
$\alpha_6, \alpha_5$ (for $\phi_3$);
whereas the automorphism \eqref{eqn:aut3} corresponds to $S_3$ symmetry
of the affine Dynkin diagram  $\Dynhat(E_6)$.  These alone do not
generate the Weyl group of $E_6$; however, together with $r_{\alpha_4}$
they do, since this extended list includes all reflections in simple roots.

\subsection{Order ideals}
\begin{definition} \rm
If $(Y, \leq)$ is a poset, an {\bf order ideal} in $Y$ is a subset 
$J \subset Y$ such that if $x \in J$, and $y \leq x$ then $y \in J$.
The set of all order ideals in $Y$ is denoted $\orderideal(Y)$ and is
itself a poset, ordered by inclusion.
\end{definition}

It is a remarkable fact that the posets $\Delta_s^+$ are related by
such a construction: there is an isomorphism 
\begin{equation}
\label{eqn:orderidealrelation}
\orderideal(\Delta_s^+) \cong
\begin{cases}
\Delta_{s'}^+ &\text{if $s = 3,4,5,6$}\\
\Delta_8^+ \setminus\{-\hat\alpha_8\} &\text{if $s = 7$}.
\end{cases}
\end{equation}
We refer the reader to \cite{Proctor} for an explanation of this
phenomenon.
Here we will explore some interesting relationships between this 
isomorphism and our compression map $f$.

\begin{definition} \rm
If $P$ and $R$ are graphs, an {\bf open map} from $P$ to $R$ is
a function $h:\text{vert}(P) \to \text{vert}(R)$ such that
\begin{enumerate}
\item $h$ is a homomorphism of graphs, i.e.  
if $(u,u') \in \text{edge}(P)$, then $(h(u), h(u')) \in \text{edge}(R)$;
\item $h$ is locally surjective, i.e. for every $v \in \text{vert}(P)$,
$h$ maps the neighbours of $v$ surjectively to the neighbours of $h(v)$.
\end{enumerate}
Equivalently, $h$ induces an open map on the topological spaces
of the graphs.
\end{definition}

\begin{proposition}
\label{prop:uniqueopenmap}
For $3 \leq s \leq 7$, there is a unique function 
$$h_s:\Delta_s^+ \to \{1, \ldots, s\}$$
such that $h_s(\alpha_s) = s$, and
$\beta \mapsto v_{h_s(\beta)}$ is an open map of graphs from 
$\Hasse_s$ to $\Dyn(E_s)$.  For $s=8$ no such function exists.
\end{proposition}

The only proof we know of this fact is to check it
case by case, which is straightforward but unenlightening.

\begin{proposition}
The isomorphism \eqref{eqn:orderidealrelation} is canonical, and
given by $\psi_s:\orderideal(\Delta_s^+) \to \Delta_{s'}^+$ where
$$
\psi_s(J) = \alpha_{s'} + \sum_{\beta \in J} \alpha_{h_s(\beta)}
$$
\end{proposition}

\begin{proof}
It is clear that any isomorphism \eqref{eqn:orderidealrelation} must
be of the form
$\psi(J) = \alpha_{s'} + \sum_{\beta \in J} \alpha_{h(\beta)}$
for some function $h:\Delta_s^+ \to \{1, \ldots, s\}$.  Since the
$\alpha_s$ is the minimal element of $\Delta_s^+$ and $\alpha_{s'}$
and $\alpha_{s'} + \alpha_s$ are the two smallest elements of
$\Delta_{s'}$, we must have $h(\alpha_s) = s$.
In light of Proposition \ref{prop:uniqueopenmap}, it suffices to
show that $h$ must induce an open map from $\Hasse_s$ to $\Dyn(E_s)$.

Suppose $\beta \succ \beta'$, is an edge of $\Hasse_s$, and let
$i = h(\beta)$, $i' = h(\beta')$.  We show that 
$(v_i, v_{i'})$ is an edge in the Dynkin diagram, i.e.
$\langle \alpha_i, \alpha_{i'}  \rangle = -1$.
Consider the order ideals
$J = \{\gamma \in \Delta_s^+\ |\ \gamma \leq \beta')$,
$J' = J \setminus \{\beta\}$ and $J'' = J \setminus \{\beta, \beta'\}$.
We have $\psi(J) \succ \psi(J') \succ \psi(J'')$, where
$\psi(J) - \alpha_i = \psi(J') 
= \psi(J'') \alpha_{i'}$.  Hence
$\langle \psi(J) - \alpha_i, \alpha_{i'} \rangle = 1$.  However, note that
$\psi(J) - \alpha_{i'}$ is not a root.  If it were then there would
two order ideals between $J$ and $J''$, namely $J'$ and 
$\psi^{-1}(\psi(J) - \alpha_{i'}$, which is impossible 
if $\beta \succ \beta'$.  Thus 
$\langle \psi(J), \alpha_{i'}  \rangle \leq 0$.  We conclude that
$\langle \alpha_i, \alpha_{i'}  \rangle
= \langle \psi(J), \alpha_{i'}  \rangle
-\langle \psi(J) - \alpha_i, \alpha_{i'}  \rangle
\leq -1$.

For local surjectivity, suppose $\beta \in \Delta_s^+$,
and let $i = h(\beta)$.  Let $(v_i, v_j)$ be an edge in the Dynkin
diagram.  We show that there exists an edge $(\beta, \gamma) \in \Hasse_s$
such that $h(\gamma) =j$.

For every $J \in \orderideal(\Delta_s^+)$,
let $J' = J \setminus \{\beta\}$, and define
$$\orderideal_\beta = \{J \in \orderideal(\Delta_s^+)
\ |\ J' \in \orderideal(\Delta_s^+)\},$$
and note that this set is non-empty.
Choose some $J \in \orderideal_\beta$.  We have 
$\psi(J) - \psi(J') = \alpha_i$, thus 
$\langle \psi(J), \alpha_j \rangle - \langle \psi(J'), \alpha_j \rangle =-1.$  
Thus either
$\langle \psi(J), \alpha_j \rangle < 0$ or
$\langle \psi(J'), \alpha_j \rangle > 0$.

In the first case, let 
$J_{\max} = \{\beta' \in \Delta_s^+\ |\ \beta' \ngtr \beta\}$ 
denote the maximal element of 
$\orderideal_\beta$.  Note that
$\langle \psi(J_{\max}), \alpha_j \rangle  \leq
\langle \psi(J), \alpha_j \rangle < 0$, and thus
$\psi(J_{\max}) + \alpha_j$ is a root.  Let 
$K = \psi^{-1}(\psi(J_{\max}) + \alpha_j)$.  We have 
$K \setminus J = \{\gamma\}$, where $h(\gamma) = j$.  Finally,
by the definition of $J_{\max}$, we must have $\gamma \succ \beta$.

In the second case, let 
$J_{\min} = \{\beta' \in \Delta_s^+\ |\ \beta' \leq \beta\}$ 
denote the maximal element of 
$\orderideal_\beta$.  We have
$\langle \psi(J_{\min}), \alpha_j \rangle \geq
\langle \psi(J), \alpha_j \rangle > 0$, and thus
$\psi(J_{\min}) - \alpha_j$ is a root.  Letting,
$K = \psi^{-1}(\psi(J_{\min}) - \alpha_j)$, we have
$J \setminus K = \{\gamma\}$, where $h(\gamma) = j$, and
$\gamma \prec \beta$.
\end{proof}

\begin{figure}[htbp]
  \begin{center}
    \epsfig{file=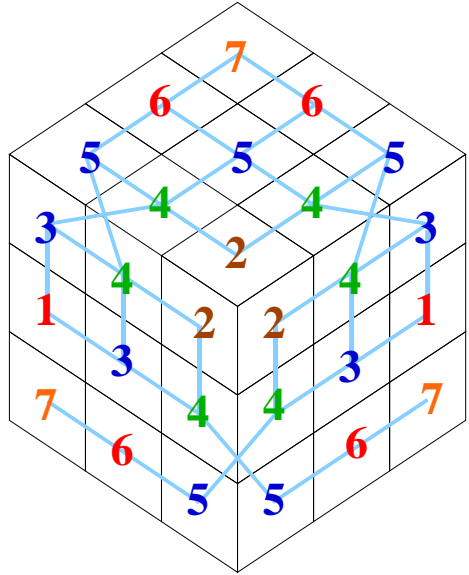, height=3in}
    \caption{The map $h_7$.}
    \label{fig:openmap7}
  \end{center}
\end{figure}

Figure \ref{fig:openmap7} shows the map $h_7$ pictured on the
corner of the $3 \times 3 \times 3$ cube.  There is a striking symmetry
in this picture.  In particular if we impose the equivalence relation
$1 \sim 6$ and $3 \sim 5$ on $\{1, \ldots, 7\}$, the numbers have
full $S_3$ symmetry.  This equivalence relation is the one that
comes from the involution on the affine Dynkin diagram $\Dynhat(E_7)$.
Furthermore the $S_3$ symmetry is exactly broken by the rule 
that $5$s and $6$s are
connected to a $7$ by a path in $\Hasse_7$ on the same face of the cube,
whereas $1$s and $3$s are not.

To understand this symmetry, let 
$\mu:\Delta_7^+ \to \Delta_7^+$, denote the automorphism from
the proof of Lemma \ref{lem:27} corresponding to $v_7 \in \Dyn(E_7)$,
and let $\rho: \Delta_7^+ \to \Delta_7^+$ denote the automorphism
corresponding to $v_6 \in \Dyn(E_7)$.
Explicitly, $\mu$ is given by the formula:
$$\mu(\beta) = -\alpha_7 - \sum_{k=1}^7 \beta^k \varepsilon(\alpha_k),$$
where $\varepsilon$ is the involution on $\Dynhat(E_7)$, taking 
$\alpha_7 \leftrightarrow \hat \alpha_7$,
$\alpha_1 \leftrightarrow \alpha_6$, $\alpha_3 \leftrightarrow \alpha_5$,
and fixing $\alpha_2$, $\alpha_4$.  
Under the identification of $\Delta_7^+$ with $\Gamma_7^+$, $\mu$
and $\rho$ are more 
simply described by $\mu(abc) = cba$, and $\rho(abc) = bca$.

The involution $\mu$ induces
an involution $\tilde \mu$ on order ideals:
$$J \leftrightarrow 
\tilde \mu(J) = \{\mu(\beta)\ |\ \beta \notin J\}.$$
Note that the root $\check \alpha_7 = f^{-1}(303) \in \Delta_7^+$, which is 
mutually orthogonal to $\alpha_7$ and $\hat \alpha_7$, is a fixed
point of $\mu$.  We use $\check \alpha_7$ to partition 
$\orderideal(\Delta_7^+)$ into four disjoint sets: 
\begin{align*}
\orderideal_0 &= \{\emptyset\} \\
\orderideal_1 &= \{J \neq \emptyset \ |\ \check \alpha_7 \notin J\} \\
\orderideal_2 &= \{J \neq \Delta_7^+ \ |\ \check\alpha_7 \in J\} \\
\orderideal_3 &= \{\Delta_7^+\}
\end{align*}
In fact $\tilde \mu$ gives a bijection between $\orderideal_i$
and $\orderideal_{3-i}$.

Note that $\orderideal_1$ can be 
identified with the order ideals in $\Gamma_7^+ \cap \antilink(033)$,
whereas $\orderideal_2$ can be identified with the order ideals
in $\Gamma_7^+ \cap \antilink(330)$.  Each order ideal in the latter
is a $120^\circ$ rotation of an order ideal in the former.
Thus we also have a bijection $\tilde \rho: \orderideal_1 \to \orderideal_2$
defined by
$$\tilde\rho(J)
= \{\rho(\beta)\ |\ \beta \in J \setminus \{\alpha_7\}\} \cup 
\{\beta\ |\ \beta \leq \check\alpha_7\}.$$

Since $\orderideal(\Delta_7^+) \cong \Delta_8^+ \setminus \{-\hat \alpha_8\}$,
we must have corresponding structures on 
$\Delta_8^+ \setminus \{-\hat \alpha_8\}$, which are related by $\psi_7$.  
Indeed we have an involution $\nu$ on this set, defined by
$$\nu(\beta) = -\beta - 2\alpha_8 - \hat \alpha_8,$$ 
where
$\psi_7 \circ \tilde \mu = \nu \circ \psi_7$.  The image of each 
$\orderideal_i$ is simply described as
$$\psi_7(\orderideal_i) = \{\beta \in \Delta_8\ |\ \beta^7 = i\}.$$
And we have a map $\sigma: \psi_7(\orderideal_1) \to \psi_7(\orderideal_2)$
given by 
$$\sigma(\beta) = -\hat\alpha_8-\mu(\beta-\alpha_8),$$ where
$\psi_7 \circ \tilde \rho = \sigma \circ \psi_7$.

We can now explain the symmetry seen in Figure~\ref{fig:openmap7}.

Suppose $J, J' \in \orderideal(\Delta_7^+)$, with 
$J = J' \cup \{\beta\}$, and let $h_7(\beta) = i$ so that 
$\psi_7(J) - \psi_7(J') = \alpha_i$.  Then
$\tilde \mu(J') = \tilde\mu(J) \cup \mu(\beta)$, and
$\psi_7(\tilde \mu(J')) - \psi_7(\tilde\mu(J)) = 
\nu(\psi_7(J')) - \nu(\psi_7(J)) = \alpha_i$.
Thus $h_7(\mu(\beta))=i= h_7(\beta)$.  This explains the reflectional
symmetry.

To see the near rotational symmetry, we
consider $J, J' \in \orderideal_1$, with $J = J' \cup \beta$. Then
$\tilde\rho(J) = \tilde\rho(J') \cup \{\rho(\beta)\}$.
We see that $\psi_7(J) - \psi_7(J') = \alpha_i$ for some $i$,
whereas 
$\psi_7(\tilde \rho(J)) - \psi_7(\tilde\rho(J')) =  
\sigma(\psi_7(J)) - \sigma(\psi_7(J')) =  
\varepsilon(\alpha_i)$.
Thus $h(\beta)$ and $h(\rho(\beta))$ will be reflections of each
other in $\Dynhat(E_7)$.

\subsection{A compression of $E_6$}
The discussion in Sections~\ref{sec:strata} and~\ref{sec:order} 
gives a description of the strata $\Delta^+_s$ for all
$s \leq 7$, but it is not the most symmetrical one for $s \leq 6$.  
For $s \leq 5$, it is easy to obtain nice
description of the strata, as they are subsets of the $D_5$
root system.  
For $s=6$, we can obtain a pleasant description by working
over $\ZZ/3$.

Let $\Delta$ be the $E_6$ root system.
Let $V = (\ZZ/3)^5$, with the standard symmetric form
\eqref{eqn:symmform3}.  Let $S = \{s_1, \ldots, s_6\}$, where
\begin{equation*}
\begin{alignedat}[b]{2}
s_1 &= 12000 &\quad s_2 &= 00012 \\
s_3 &= 01200 &\quad s_4 &= 00120 \\
s_5 &= 00011 &\quad s_6 &= 11111 
\end{alignedat}
\qquad \qquad
\epsfig{file=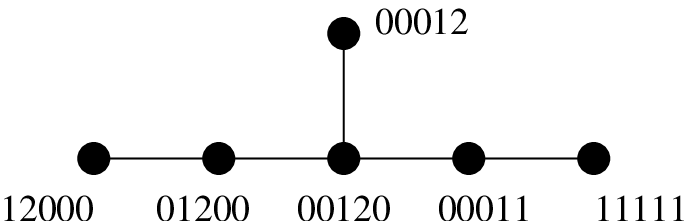, height=0.6in}
\end{equation*}

\begin{proposition}
For $S$ as above the relations \eqref{eqn:samedynkin} hold.
\end{proposition}

As a consequence of we obtain the following corollary of Theorem
\ref{thm:injective}, we have

\begin{corollary}
The map $f: \Delta \to \Gamma$ is a bijection.
\end{corollary}

\begin{proof}
It is an injection by Theorem \ref{thm:injective}.  To show it
is a bijection, we must calculate the size of $\Gamma$.  The
vectors in $\Gamma$ have either $2$ or $5$ non-zero coordinates,
each of which can be $\pm 1$.
Thus there are $2^2 {5 \choose 2} + 2^5 {5 \choose 5} = 72$ 
elements in $\Gamma$.  But $\#(\Delta) = 72$, so $f$ is a 
bijection.
\end{proof}

Let $\Gamma_s = f(\Delta_s)$, and $\Gamma_s^+ = f(\Delta_s^+)$
denote the images of the strata under $f$.

\begin{theorem}
The image of the top stratum, $\Gamma^+_6$, is the set of vectors in 
$V$ with all 
coordinates non-zero, and an even number of coordinates equal to $2$:
$$
\Gamma^+_6 = 
\Big\{(x_1, \ldots, x_5) \in V\ \Big|\ \prod_{i=1}^5 x_i = 1\Big\}$$
\end{theorem}

\begin{proof}

The argument is parallel to the proof of Lemma \ref{lem:27}.
One can check that the automorphisms 
corresponding to Dynkin diagram vertices $v_3, v_4, v_5$ give permutations 
of the coordinates which generate the symmetric group $S_5$,
hence 
$S_5$ acts on the coordinates of $\Gamma_6^+$.  Furthermore, the
automorphism corresponding to $v_6$ is $(x_1, x_2, x_3, x_4, x_5)
\mapsto (-x_1, -x_2, -x_3, -x_4, x_5)$.  Applying these automorphisms
to $11111 = f(\alpha_6)$, we see that all 16 elements of $\Gamma_6^+$
are indeed of the correct form.
\end{proof}

\begin{definition} \rm
The {\bf T-graph} for $V$ is the graph whose vertex set is $V$,
and $x= (x_1, \ldots, x_5)$ is 
adjacent to $y= (y_1, \ldots, y_5)$
if $x_i = y_i$ for exactly one $i$.
\end{definition}

\begin{theorem}
\label{thm:T-graph6}
The T-graph and the O-graph agree when restricted to the image of
any stratum $f(\Delta_s^+)$.
\end{theorem}

The proof is analogous to that of Theorem \ref{thm:T-graph7}.

\begin{theorem}
\label{thm:order6}
Let $\alpha$ be a positive root.  Let 
$\beta \in \Delta_s^+$.  Then
$f(\beta)+ f(\alpha) \in \Gamma_s^+$ if and only if
$\beta + \alpha \in \Delta_s^+$.
\end{theorem}

\begin{proof}
Certainly if  $\beta + \alpha \in \Delta_s^+$, then
$f(\beta) + f(\alpha) \in \Gamma_s^+$.  Suppose that
$\beta + \alpha \notin \Delta_s^+$.  If $\beta+\alpha$
is a root, then it belongs to
some stratum other than $\Delta_s^+$, so 
$f(\beta) + f(\alpha) \notin \Gamma_s^+$.
If $\beta = \alpha$, then 
$f(\beta) + f(\alpha) = f(-\beta)$ is the image of
a negative root, so $f(\beta) + f(\alpha) \notin \Gamma_s^+$.
Otherwise, since $\beta + \alpha$ is not a root, we
must have
$\langle \beta, \alpha \rangle \in \{0,1\}$.
In this case, 
$\langle \beta+\alpha, \beta+\alpha \rangle 
= 4 + 2\langle \beta, \alpha\rangle \neq 2
\pmod 3$, so in fact
$f(\beta) + f(\alpha) \notin \Gamma$.
\end{proof}

\begin{figure}[htbp]
  \begin{center}
    \epsfig{file=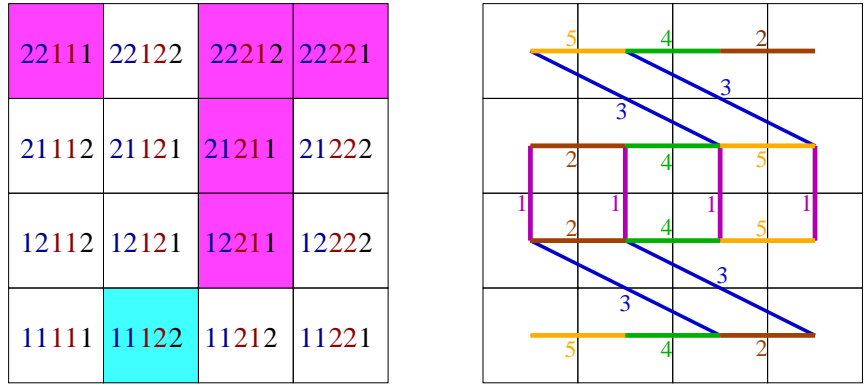, height=2in}
    \caption{Orthogonality and partial order in $\Gamma_6^+$. }
    \label{fig:square}
  \end{center}
\end{figure}

As we did with $\Delta_7^+$, we can recover the partial order
structure on $\Delta_6^+$ (and the smaller strata) using 
Theorem \ref{thm:order6}.  If the elements of $\Gamma_6^+$
are arranged as shown in Figure~\ref{fig:square}, we join
$x$ and $y$ if $x - y$ is a simple root.  In light of
Theorem \ref{thm:T-graph6}, orthogonality is also easily
determined in this picture.  Figure~\ref{fig:square} shows
the example of $\link(11122) \cap \Gamma_6^+$, which
gives the set of root images orthogonal to $11122$.

\end{document}